\documentclass[10pt]{article}
\usepackage{amsmath,amsthm,amsfonts,amssymb}
\usepackage{epsfig}
\usepackage{color}
\usepackage{latexsym}
\usepackage{supertabular}
\usepackage{graphicx}
\usepackage{a4wide}
\numberwithin{equation}{section}

\def\whitebox{{\hbox{\hskip 1pt
 \vrule height 6pt depth 1.5pt
 \lower 1.5pt\vbox to 7.5pt{\hrule width
    3.2pt\vfill\hrule width 3.2pt}%
 \vrule height 6pt depth 1.5pt
 \hskip 1pt } }}
\def\qed{\ifhmode\allowbreak\else\nobreak\fi\hfill\quad\nobreak
     \whitebox\medbreak}

\newcommand{\ignore}[1]{}

\theoremstyle{plain}
\newtheorem{theorem}{Theorem}[section]
\newtheorem{corollary}[theorem]{Corollary}
\newtheorem{lemma}[theorem]{Lemma}

\newtheorem{example}[theorem]{Example}

\newtheorem{remark}[theorem]{Remark}
\newtheorem{conjecture}[theorem]{Conjecture}
\def\qed{{\hfill$\square$}}
\def\proof{{\vspace{-0.3cm}\bf Proof: \,}}
\def\N{{\mathbb N}}
\def\Z{{\mathbb Z}}
\def\Q{{\mathbb Q}}

\def\C{{\mathbb C}}
\def\F{{\mathbb F}}
\def\mod{{\mathrm{mod\,\,}}}

\def\Gal{{\mathrm{Gal}}}

\def\Tr{{\mathrm{Tr}}}
\def\Norm{{\mathrm{Norm}}}
\def\Cay{{\mathrm{Cay}}}
\def\ord{{\mathrm{ord}}}

\title{Strongly Regular Cayley Graphs, Skew Hadamard Difference Sets, 
and Rationality of Relative Gauss Sums}
\author{Koji Momihara\footnotemark[1]}
\date{} 
\begin{document}
\maketitle
\footnotetext[1]{
Department of Mathematics, Faculty of Education, Kumamoto University,  
2-40-1 Kurokami, Kumamoto 860-8555, Japan; Email address: 
momihara@educ.kumamoto-u.ac.jp}
\renewcommand{\thefootnote}{\arabic{footnote}}
\begin{abstract}  
In this paper, we give constructions of strongly regular Cayley graphs and skew Hadamard difference sets. 
Both constructions are based on choosing  cyclotomic classes in finite fields, and our results generalize ten of the eleven sporadic examples of cyclotomic strongly regular graphs given by  Schmidt and White \cite{SW02} and several of subfield examples into infinite families. 
These infinite families of strongly regular graphs have new parameters. The main tools that we employed are relative Gauss sums instead of explicit evaluations of Gauss sums. 
\end{abstract}
\begin{center} 
\end{center}
\section{Introduction}
In this paper, we will assume that the reader is familiar with the theory of strongly regular graphs and difference sets. For the theory of strongly regular graphs (srgs), 
our main reference is  the lecture note of Brouwer and Haemers \cite{bh}. For difference sets, we refer the reader to 
Chapter 6 of \cite{bjl}. 
We remark that strongly regular graphs are closely related to other combinatorial objects, such as two-weight codes, two-intersection sets in finite geometry, 
and partial difference sets. For these connections, we refer the reader to \cite[p.~132]{bh}, \cite{CK86, M94}.  

Let $\Gamma$ be a simple and undirected graph and $A$ be its adjacency matrix. A useful way to check whether a graph is strongly regular is by using the eigenvalues of its adjacency matrix. For convenience we call an eigenvalue {\it restricted} if it has an eigenvector perpendicular to the all-ones vector ${\bf 1}$. (For a $k$-regular connected graph, the restricted eigenvalues are the eigenvalues different from $k$.)

\begin{theorem}\label{char}
For a simple graph $\Gamma$ of order $v$, not complete or edgeless, with adjacency matrix $A$, the following are equivalent:
\begin{enumerate}
\item $\Gamma$ is strongly regular with parameters $(v, k, \lambda, \mu)$ for certain integers $k, \lambda, \mu$,
\item $A^2 =(\lambda-\mu)A+(k-\mu) I+\mu J$ for certain real numbers $k,\lambda, \mu$, where $I, J$ are the identity matrix and the all-ones matrix, respectively, 
\item $A$ has precisely two distinct restricted eigenvalues.
\end{enumerate}
\end{theorem}
One of the most effective methods for constructing srgs is by the Cayley graph construction. For example, the Paley graph ${\rm P}(q)$ is one class of well known Cayley graphs, that is, the graph with the finite field $\F_q$ as vertex set, where two vertices are adjacent when 
their difference is a nonzero quadratic. It has the parameters 
$(v,k,\lambda,\mu)=(4t+1,2t,t-1,t)$. 
In general, 
let $G$ be an additively written group of order $v$, and let $D$ be a subset of $G$ such that $0\not\in D$ and $-D=D$, where $-D=\{-d\mid d\in D\}$. The {\it Cayley graph on $G$ with connection set $D$}, denoted ${\rm Cay}(G,D)$, is the graph with the elements of $G$ as vertices; two vertices are adjacent if and only if their difference belongs to $D$. In the case when $\Cay(G,D)$ is strongly regular, the connection set $D$ is called a (regular) {\it partial difference set}.  The survey of Ma~\cite{M94} contains much of what is known about partial difference sets and about connections with strongly regular Cayley graphs. 

A difference set $D$ in an (additively written) finite group $G$ is called {\it skew Hadamard} if $G$ is the disjoint union of $D$, $-D$, and $\{0\}$. The primary example (and for many years, the only known example in abelian groups) of skew Hadamard difference sets is the classical Paley difference set in $(\F_q,+)$ consisting of the nonzero squares of $\F_q$, where $\F_q$ is the finite field of order $q$, a prime power congruent to 3 modulo 4.  Skew Hadamard difference sets are currently under intensive study; see the introduction of \cite{FX113} for a short survey of known constructions of skew Hadamard difference sets and related problems.
As we see in the next section, in order to check that  a candidate subset $D$ of $\F_q$ is a partial difference set or a skew Hadamard
difference set in $(\F_q, +)$, it is sufficient to compute certain  character sums of $\F_q$ in common. 

A classical method for constructing both connection sets of strongly regular graphs (i.e., partial difference sets) and difference sets in the additive groups of finite fields is to use cyclotomic classes of finite fields. 
Let $p$ be a prime,  $f$ a positive integer, and let $q=p^f$. Let $k>1$ be an integer such that $k|(q-1)$, and $\gamma$ be a primitive element of $\F_q$. 
Then the cosets $C_i^{(k,q)}=\gamma^i \langle \gamma^k\rangle$, $0\leq i\leq k-1$, are called the {\it cyclotomic classes of order $k$} of $\F_q$. 
Many authors have studied the problem of determining when a union $D$ of some cyclotomic classes forms a (partial) difference set.
Especially, when $D$ consists of only a subgroup of $\F_q$, 
many authors have studied extensively \cite{BM73,BMW82,bwx,FX111,FX113,FMX11,GXY11,IM10,DL95,
VLSch,M75,SW02,S67}. 
(Some of these authors used the language of cyclic codes
in their investigations instead of strongly regular Cayley graphs or partial difference sets. We choose to use the language of srg.)
We call a strongly regular Cayley graph $\Cay(\F_q,D)$ {\it cyclotomic} if 
$D$ is such. 
The Paley graphs are primary examples of  cyclotomic srgs.  
Also, if $D$ is the multiplicative group of a
subfield of $\F_q$, then it is clear that $\Cay(\F_q , D)$ is strongly regular. 
These cyclotomic srgs are usually called {\it subfield examples}. Next, if there exists a
positive integer $t$ such that $p^t\equiv -1\,(\mod{k})$, then $\Cay(\F_q , D)$ is strongly regular. This case had already generalized so that 
$D$ is a union of some cyclotomic cosets based on the computation of ``pure Gauss sums", see \cite{bwx,VLSch}. 
These examples are usually called {\it semi-primitive}. 
Schmidt and White presented the following conjecture on cyclotomic srgs. 
\begin{conjecture}\label{con:SW}(\cite{SW02})
Let $\F_{p^f}$ be the finite field, $k\,|\,\frac{p^f-1}{p-1}$ with $k>1$, and $C_0:=C_{0}^{(k,p^f)}$ with $-C_0=C_0$. If $\Cay(\F_{p^f},C_0)$ is strongly regular, then one of the following holds: 
\begin{enumerate}
\item[(1)] (subfield case) $C_0=\F_{p^d}^\ast$ where $d\,|\,f$,  
\item[(2)] (semi-primitive case) $-1\in \langle p\rangle\le (\Z/k\Z)^\ast$,
\item[(3)] (exceptional case) $\Cay(\F_{p^f},C_0)$ has one of the parameters 
given in Table~\ref{Tab1}. 
\begin{table}[h]
\caption{Eleven sporadic examples}
\label{Tab1}
$$
\begin{array}{|c||c|c|c|c|}
\hline
\mbox{No.}&k&p&f&e:=[(\Z/k\Z)^\ast:\langle p\rangle]\\
\hline
1&11&3&5&2\\
2&19&5&9&2\\
3&35&3&12&2\\
4&37&7&9&4\\
5&43&11&7&6\\
6&67&17&33&2\\
7&107&3&53&2\\
8&133&5&18&6\\
9&163&41&81&2\\
10&323&3&144&2\\
11&499&5&249&2\\
\hline
\end{array}
$$
\end{table}
\end{enumerate}
\end{conjecture}
Recently, the authors of \cite{FX113,FMX11,GXY11} succeeded to generalize the examples of Table~\ref{Tab1} except for srgs of No. 1, 5,  and 8 into 
infinite families using ``index $2$ and $4$ Gauss sums''. 
\begin{theorem} \label{pro:knows}
\begin{itemize}
\item[(i)] (\cite{FX111}) Let $q=p^{p_1^{m-1}(p_1-1)/2}$, $k=p_1^m$, and  $D=\bigcup_{i=0}^{p_1^{m-1}-1}C_{i}^{(k,q)}$. Then, $\Cay(\F_q,D)$ is strongly regular for any $m$ in 
the following cases:  
\begin{eqnarray*}
(p,p_1)=(2,7),(3,107),(5,19),(5,499),(17,67),(41,163). 
\end{eqnarray*}
\item[(ii)] (\cite{GXY11}) Let $q=p^{p_1^{m-1}(p_1-1)/4}$, $k=p_1^m$, and  $D=\bigcup_{i=0}^{p_1^{m-1}-1}C_{i}^{(k,q)}$. 
Then, $\Cay(\F_q,D)$ is strongly regular for any $m$ in 
the following cases:  
\begin{eqnarray*}
(p,p_1)=(3,13),(7,37). 
\end{eqnarray*}
\item[(iii)] (\cite{FMX11}) Let $q=p^{p_1^{m-1}(p_1-1)p_2^{n-1}(p_2-1)/2}$, $k=p_1^mp_2^n$, and  
$D=\bigcup_{i=0}^{p_1^{m-1}-1}\bigcup_{j=0}^{p_2^{n-1}-1}
C_{p_2^n i+p_1^m j}^{(k,q)}$. Then, $\Cay(\F_q,D)$ is strongly regular for any $n$ and $m$ in 
the following cases:  
\begin{eqnarray*}
(p,p_1,p_2)=(2,3,5),(3,5,7),(3,17,19). 
\end{eqnarray*}
\end{itemize}
\end{theorem}
The srgs in the cases when $(p,p_1)=(2,7)$ of (i), $(p,p_1)=(3,13)$ of (ii), and $(p,p_1,p_2)=(2,3,5)$ of (iii) of Theorem~\ref{pro:knows} are 
generalizations of subfield examples. The others are generalizations of  sporadic examples of Table~\ref{Tab1}.
Note that it is impossible to generalize the example of No. 1 of Table~\ref{Tab1} by a similar manner since  $\langle 3\rangle\le (\Z/11^m\Z)^\ast$ is not of index $2$ for $m>1$. 

In \cite{FX113,FMX11}, the following two constructions of skew Hadamard  difference sets 
and Paley type partial difference sets were given. 
(A partial difference set $D$ in a group $G$ is said to be of 
{\it Paley type} if the parameters of the corresponding strongly 
regular Cayley graph are $(v, \frac{v-1}{2}, \frac{v-5}{4}, \frac{v-1}{4})$.) 
\begin{theorem}\label{IntroThm3}
\begin{itemize}
\item[(i)] (\cite{FX113}) Let $p_1\equiv 7\,(\mod{8})$ be a prime, $k=2p_1^m$, and 
let $p$ be a prime such that $f:=\ord_k(p)=\phi(k)/2$, where $\phi$ is the Euler totient function. Let $s$ be an odd integer, $H$ denotes any subset of $\Z_k$ such that 
$\{i\,(\mod{p_1^m})\,|\,i\in H\}=\Z_{p_1^m}$, and  let $D=\bigcup_{i\in H}C_i^{(k,p^{fs})}$. Then, $D$ is a skew Hadamard 
difference set if $p\equiv 3\,(\mod{4})$ and $D$ is a 
Paley type partial difference set if $p\equiv 1\,(\mod{4})$. 
\item[(ii)]  (\cite{FMX11})  Let $q=p^{p_1^{m-1}(p_1-1)/2}$, $k=2p_1^m$, and $H=Q\cup 2Q\cup \{0\}$, where 
$Q$ is the subgroup of index $2$ of $(\Z/2p_1\Z)^\ast$. Set 
$D=\bigcup_{j=0}^{p_1^{m-1}}\bigcup_{i\in H}C_{2j+ip_1^{m-1}}^{(k,q)}$. 
Then, $D$ is a skew Hadamard difference set in the case when $(p,p_1)=(3,107)$ and $D$ is  
a Paley type partial difference set in the  cases when 
\[
(p,p_1)=(5,19),(17,67),(41,163),(5,499). 
\]
\end{itemize}
\end{theorem}
The proofs of the above theorems are based on index $2$ and $4$ Gauss sums. In order to show that the srgs of No. 5 and 8 in Table~\ref{Tab1} lead to infinite families, we need to explicitly evaluate index $6$ Gauss sums if we apply a similar technique of  
\cite{FX111,FX113,FMX11,GXY11}. 
However, it seems to be difficult to compute index more than $4$ Gauss sums, and this implies that it is hard to find new strongly regular graphs or skew Hadamard difference sets on $\F_q$ from 
index more than $4$ cases.  
In this paper, we will show that  explicit evaluations of 
Gauss sums are not needed if some initial examples of strongly regular Cayley graphs or skew Hadamard difference sets  satisfying  certain conditions are 
found. Instead, 
we will investigate the rationality of ``relative Gauss sums". 
As consequences, we generalize the srgs of No. 5 and 8 in Table~\ref{Tab1} into infinite families and find further infinite families 
of cyclotomic srgs with new parameters as generalizations of 
subfield examples (see Tables~\ref{Tab4} and \ref{Tab6} in Section~\ref{sec:stro}). Furthermore, 
we obtain 
two infinite families of skew Hadamard difference sets in $(\F_q,+)$, where $q=3^{3\cdot 13^{m-1}}$ and $7^{7\cdot 29^{m-1}}$. 
\section{Rationality of relative Gauss sums}\label{RGauss}
\subsection{Preliminary}
Let $p$ be a prime, $f$ a positive integer, and $q=p^f$. The canonical additive character $\psi$ of $\F_q$ is defined by 
$$\psi\colon\F_q\to \C^{\ast},\qquad\psi(x)=\zeta_p^{\Tr _{q/p}(x)},$$
where $\zeta_p={\rm exp}(\frac {2\pi i}{p})$ and $\Tr _{q/p}$ is the trace from $\F_q$ to $\F_p$. For a multiplicative character 
$\chi$ of $\F_q$, we define the {\it Gauss sum} 
\[
G_f(\chi)=\sum_{x\in \F_q^\ast}\chi(x)\psi(x), 
\] 
which belongs to $\Z[\zeta_{kp}]$ of integers in the cyclotomic field $\Q(\zeta_{kp})$, where $m$ is the order of $\chi$.
Let $\sigma_{a,b}$ be the automorphism of $\Q(\zeta_{kp})$ determined 
by 
\[
\sigma_{a,b}(\zeta_k)=\zeta_{k}^a, \qquad
\sigma_{a,b}(\zeta_p)=\zeta_{p}^b 
\]
for $\gcd{(a,k)}=\gcd{(b,p)}=1$. 
Below are several basic properties of Gauss sums \cite{LN97}: 
\begin{itemize}
\item[(i)] $G_f(\chi)\overline{G_f(\chi)}=q$ if $\chi$ is nontrivial;
\item[(ii)] $G_f(\chi^p)=G_f(\chi)$, where $p$ is the characteristic of $\F_q$; 
\item[(iii)] $G_f(\chi^{-1})=\chi(-1)\overline{G_f(\chi)}$;
\item[(iv)] $G_f(\chi)=-1$ if $\chi$ is trivial;
\item[(v)] $\sigma_{a,b}(G_f(\chi))=\chi^{-a}(b)G_f(\chi^a)$;
\end{itemize}

In general, the explicit evaluation of Gauss sums is a very difficult problem. There are only a few cases where the Gauss sums have been evaluated. 
The most well known case is {\it quadratic} case, in other words, the order of $\chi$ is two. 
\begin{lemma}(\cite{LN97})\label{le:quad}
Let $\eta$ be the quadratic character of $\F_q=\F_{p^f}$. Then, it holds that 
\[
G_f(\eta)=(-1)^{f-1}\left(\sqrt{(-1)^{\frac{p-1}{2}}p}\right)^f. 
\]
\end{lemma}
The next simple case is the so-called {\it semi-primitive case} (also 
referred to as {\it uniform cyclotomy} or {\it pure Gauss sum}), where there 
exists an integer $j$ such that $p^j\equiv -1\,(\mod{k})$, where $k$ is the order of 
the multiplicative character $\chi$ involved. 
\begin{theorem}\label{thm:semiprim}(\cite{BEW97})
Suppose that $k>2$ and $p$ is semi-primitive modulo $k$, 
i.e., there exists an $s$ s.t. $p^s\equiv -1\,(\mod{k})$. Choose 
$s$ minimal and write 
$f=2st$. Let $\chi$ be a multiplicative character of order $k$. 
Then, 
\[
p^{-f/2}G_f(\chi)=
\left\{
\begin{array}{ll}
(-1)^{t-1}&  \mbox{if $p=2$;}\\
(-1)^{t-1+(p^s+1)t/k}&  \mbox{if $p>2$. }
 \end{array}
\right.
\]
\end{theorem}
This theorem was used to find strongly regular graphs and difference sets 
on $\F_q$, e.g., see \cite{BMW82,bwx}. 

The next interesting case is the index $2$ case where the subgroup $\langle p\rangle$ generated by $p\in (\Z/{k}\Z)^\ast$ has index $2$ in $(\Z/{k}\Z)^\ast$ and $-1\not\in \langle p\rangle $. In this case, 
it is known that $k$ can have at most two odd prime divisors. 
Many authors have investigated this case, see e.g., \cite{BM73,L97,M98,MV03,YX10,YX11}. 
In particular, complete solution to the problem of evaluating Gauss sums in this case was 
recently given in \cite{YX10}. Also, the index $4$ case was treated in 
\cite{FYL05,YLF06}. Recently, these index $2$ and $4$ Gauss sums were applied to 
show the existence of infinite families of 
new strongly regular graphs and skew Hadamard difference sets 
on $\F_q$ in \cite{FX111,FX113,FMX11,GXY11}. 
However, it is quite difficult to  
explicitly evaluate Gauss sums of general index. This 
implies that it is difficult to find new strongly regular graphs on $\F_q$ from 
index more than $4$ cases if we apply a similar technique of  
\cite{FX111,FX113,FMX11,GXY11}. However, we will show in Section~\ref{sec:const} of this paper that explicit evaluations of 
Gauss sums are not needed if some initial examples of strongly regular graphs or skew Hadamard difference sets  satisfying  certain conditions are 
found. Instead, 
we will use rationality of relative Gauss sums. 
For two nontrivial multiplicative characters $\chi$ of $\F_{p^f}$ and $\chi'$ of $\F_{p^{f'}}$ with $f\,|\,f'$, the {\it relative Gauss sum associated with $\chi$ and $\chi'$} is defined 
as
\[
\vartheta_p(\chi',\chi):=\frac{G_{f'}(\chi')}{p^{\frac{f'-f}{2}}G_f(\chi)}. 
\] 
In particular, we investigate when $\vartheta_p(\chi',\chi)=1$ or $-1$ holds in the case where both of 
$G_{f'}(\chi')$ and $G_f(\chi)$ are of index $e$ case. 
Note that the concept of relative Gauss sums was introduced in \cite{Y85} 
as the fractional $G_{f'}(\chi')/G_f({\chi})$, where $\chi$ is the restriction of 
$\chi'$ to $\F_{p^f}$. Hence, our definition generalize his definition and 
normalize so that the absolute value is equal to $1$ when $\chi$ and $\chi'$ are nontrivial. 

Below, we give important formulae on Gauss sums. 
The following is known as  {\it the Davenport-Hasse lifting formula}.
\begin{theorem}\label{thm:lift}(\cite{BEW97,LN97})
Let $\chi$ be a nontrivial character on $\F_q=\F_{p^f}$ and 
let $\chi'$ be the lifted character of $\chi$ to $\F_{q'}=\F_{p^{fs}}$, i.e., $\chi'(\alpha):=\chi(\Norm_{\F_{q'}/\F_q}(\alpha))$ for $\alpha\in \F_{q'}$. 
Then, it holds that 
\[
G_{fs}(\chi')=(-1)^{s-1}(G_f(\chi))^s. 
\]
\end{theorem}
The following is  called {\it the Davenport-Hasse product 
formula}. 
\begin{theorem}\label{thm:Stickel2}(\cite{BEW97})
Let $\eta$ be a character on $\F_q=\F_{p^r}$ of order $\ell>1$. For  every nontrivial character $\chi$ on $\F_q$, 
\[
G_r(\chi)=\frac{G_r(\chi^\ell)}{\chi^\ell(\ell)}
\prod_{i=1}^{\ell-1}
\frac{G_r(\chi\eta^i)}{G_r(\eta^i)}. 
\]
\end{theorem}
We close this subsection providing the following lemma~\cite{Y85}.  
\begin{lemma}\label{Prop:Yamamoto}(\cite{Y85})
Let $\chi'$ be a character of order $k'$ of $\F_{p^{f'}}$ and $\chi$ be the restriction of $\chi' $ to $\F_{p^f}$, where $f\,|\,f'$. If $\chi$ is nontrivial on 
$\F_{p^f}$, it holds that  
\[
p^\frac{f-f'}{2}\vartheta_p(\chi',\chi)=\sum_{x\in L; \Tr_{p^{f'}/p^f}(x)=1}\chi'(x),  
\]
where $L$ is a set of representatives for $\F_{p^{f'}}^\ast/\F_{p^{f}}^\ast$. 
\end{lemma}
\subsection{Relative Gauss sums}
In this section, fix an integer $k>1$, and let  $p$ be a prime such that 
$\gcd{(p,k)}=1$. Let $f$ be the order of $p$ in $(\Z/k\Z)^\ast$ and  set 
$q=p^f$.
Write $\zeta_k=e^{2\pi i/k}$ and $\zeta_p=e^{2\pi i/p}$. 
Define 
\[
K=\Q(\zeta_k),\, M=K(\zeta_p)=\Q(\zeta_k,\zeta_p), 
\]
and 
let $O_k$ and $O_M$ denote their respective rings of integers. 
For $j\in (\Z/k\Z)^\ast$, define $\sigma_j\in \Gal(M/\Q(\zeta_p))$ by 
$\sigma_j(\zeta_k)=\zeta_k^j$. Let $P$ be a prime ideal of $O_K$ lying over 
$p$. Then, for some prime ideal ${\mathfrak{p}}$ of $O_M$ such that 
$PO_M={\mathfrak{p}}^{p-1}$ and ${\mathfrak{p}}\cap O_K=P$. 
Write $P_j=\sigma_j(P)$ and ${\mathfrak{p}}_j=\sigma_j({\mathfrak{p}})$,  and 
then $P_jO_M={\mathfrak{p}}_j^{p-1}$. Let $T$ be a set of representatives 
of
$(\Z/{k}\Z)^\ast/\langle p\rangle$. Then, $pO_K=\prod_{j\in T}P_j$ follows, 
where 
$P_j$ are all distinct, and hence $pO_M=\prod_{j\in T}{\mathfrak{p}}_j^{p-1}$ holds.

Define the character $\chi_P$ of order $k$ on the finite field 
$O_K/P$ by letting $\chi_P(\alpha+P)$ denote the 
unique power of $\zeta_k$ such that 
\[
\chi_P(\alpha+P)\equiv \alpha^{(q-1)/k}\, (\mod{P}), 
\] 
when $\alpha\in O_K\setminus P$. When $\alpha\in P$, set $\chi_P(\alpha+P)=0$. We call $\chi_P$ {\it the 
Teichm\"{u}ller character 
associated to $P$}. 
Now we identify $\chi_P$ with a character of 
$\F_q$. 

Define 
\[
\theta(k,p)=\sum_{t\in (\Z/{k}\Z)^\ast}\left\langle \frac{t}{k}\right\rangle
\sigma_t^{-1}, 
\]
called {\it the Stickelberger element}, where 
$\left\langle x\right\rangle$ is the fractional part of the rational $x$. 
 Every 
integer $a$ can be written uniquely in the form $\sum_{i=0}^na_ip^i$, where 
$0\le a_i<p$. We denote by $s_p(a)$ the sum of all $a_i$. The following 
are given in \cite{BEW97,IR}. 
\begin{lemma}
For any integer $a$, $0\le a<q-1$, we have 
\[
s_p(a)=(p-1)\sum_{i=0}^{f-1}\left\langle \frac{p^i a}{q-1}\right\rangle.
\] 
\end{lemma}  
\begin{theorem}
Let $k$ be a positive integer. Let $p$ be a prime such that $\gcd{(p,k)}=1$  and $f$ be the order of $p$ 
in $(\Z/{k}\Z)^\ast$. For a prime ideal  
${\mathfrak{p}}$ of $O_M$ lying over $P$, it holds 
\[
G_f(\chi_P^{-1})O_M={\mathfrak{p}}^{\sum_{t\in T}s_p(t(q-1)/k)\sigma_t^{-1}}
={\mathfrak{p}}^{(p-1)\sum_{t\in T}\sum_{i=0}^{f-1}
\langle tp^i/k\rangle\sigma_t^{-1}}\subseteq O_{M}. 
\] 
\end{theorem}
This theorem is known as {\it the Stickelberger relation}. By the relation 
$G_f(\chi^a)G_f(\chi^{-a})=\pm p^f$, we also have 
\[
G_f(\chi_P)O_M
={\mathfrak{p}}^{(p-1)(f\sum_{t\in T}\sigma_t-\sum_{t\in T}\sum_{i=0}^{f-1}
\langle tp^i/k\rangle\sigma_t^{-1})}. 
\]

In the rest of this paper, we will assume the following. 
Let $h=2^tp_1p_2\cdots p_\ell$ be a positive integer with distinct odd primes  $p_i$ and $p$ be a prime satisfying the following:
For any divisor $d=2^s p_{i_1}\cdots p_{i_m}$ of $h$, if 
$\langle p\rangle$ is of index $u$ modulo $d$, then so does $\langle p\rangle$ modulo $d'=2^s p_{i_1}^{x_1}\cdots p_{i_m}^{x_m}$ for any $x_i\ge 1$. 
Let $e$ denotes the index of $\langle p\rangle$ modulo $h$. 
Let $p_1$ be an odd prime factor of $k=2^tp_1^{e_1}p_2^{e_2}\cdots p_\ell^{e_\ell}$ and 
set $k'=kp_1$. Then, by the assumption, $\langle p\rangle$ is again of index $e$ in both of  
$(\Z/k\Z)^\ast$ and $(\Z/k'\Z)^\ast$. 
Set $q=p^{f}$ and $q'=p^{f'}$, where $f=\phi(k)/e$ and 
 $f'=\phi(k')/e$. 

Let $O_K$, $O_{K'}$, $O_M$, $O_{M'}$, $O_L$, $O_{L'}$ denote the respective rings of integers 
of $\Q(\zeta_k)$, $\Q(\zeta_{k'})$, $\Q(\zeta_k,\zeta_p)$, $\Q(\zeta_{k'},\zeta_{p})$, $\Q(\zeta_{p^{f}-1})$, $\Q(\zeta_{p^{f'}-1})$. 
Let 
$P\subseteq O_K$ be a prime ideal lying over $p$ and 
${\mathfrak{p}}\subseteq O_M$ be a prime ideal lying over $P$. Also,  let ${\mathfrak{p}}'\subseteq O_{M'}$ be a prime ideal lying over  
${\mathfrak{p}}$ and let $P'={\mathfrak{p}}'\cap O_{K'}$, so that 
${\mathfrak{p}}'\cap O_{M}={\mathfrak{p}}$ and 
$P'\cap O_K=P$. 
Let $T'$ be a set of representatives for $(\Z/k'\Z)^\ast/\langle p\rangle$. Then, there is a one to one correspondence between 
$\{\sigma_j(P)(=:P_j)\,|\,j\in T\}$ and $\{\sigma_j'(P)(=:P_j')\,|\,j\in T'\}$ such that $P_j=P_j'\cap O_K$, where $\sigma_j'\in \Gal(\Q(\zeta_{k'p})/\Q(\zeta_p))$ satisfying $\sigma_j'(\zeta_{k'})=\zeta_{k'}^{j}$.  
By multiplying $O_{K'}$ to both side of 
$pO_K=\prod_{j\in T}P_j$, together with $pO_{K'}=\prod_{j\in T'}P_j'$,  we have 
$P_jO_{K'}=P_j'$. Furthermore, by multiplying 
$O_{M'}$ to both side of $P_jO_{K'}=P_j'$, we have 
$P_jO_{M'}=P_j'O_{M'}={\mathfrak{p'}}_j^{p-1}$, where ${\mathfrak{p'}}_j\subseteq O_{M'}$ is a prime ideal lying over $P_j'$. On the other hand, since 
$P_jO_{M'}={\mathfrak{p}}_j^{p-1} O_{M'}$, we obtain 
${\mathfrak{p}}_jO_{M'}={\mathfrak{p'}}_j$. 

Let ${\mathfrak{P}}\subseteq O_{L}$ and ${\mathfrak{P'}}\subseteq O_{L'}$
be prime ideals lying over $P$ and $P'$, respectively. 
It is known that $O_{L}/{\mathfrak{P}}=\{\alpha+{\mathfrak{P}}\,|\,\alpha\in O_K/P\}$ and that 
\[
\chi_{{\mathfrak{P}}}^\frac{p^{f}-1}{k}(\alpha+{\mathfrak{P}})=\chi_P(\alpha+P)
\]
for $\alpha\in O_K$, so that 
\[
G_f(\chi_{{\mathfrak{P}}}^{a\frac{p^f-1}{k}})=G_f(\chi_P^a).
\]
See Exercise~11-1 of \cite{BEW97}. 
Now, we can take the set $\{0\}\cup 
\{\zeta_{p^{f}-1}^i\,|\,0\le i\le p^{f}-1\}$ as representatives for $O_{L}/{\mathfrak{P}}$ and then 
\[
\chi_{{\mathfrak{P}}}(\zeta_{ p^{f}-1}^i+{\mathfrak{P}})=\zeta_{p^{f}-1}^i.
\]  

By the definition of Teichm{\"u}ller characters, for $\alpha\in (\zeta_{p^{f}-1}+{\mathfrak{P}})\cap O_{K}$ and $\beta\in 
(\zeta_{p^{f'}-1}+{\mathfrak{P'}})\cap O_{K'}$ it holds 
\begin{eqnarray}\label{eq:basicchar}
\chi_P(\alpha^i+P)&=&\chi_{\mathfrak{P}}(\zeta_{p^{f}-1}^{\frac{p^f-1}{k}i}+{\mathfrak{P}})
=\chi_{\mathfrak{P'}}^{\frac{p^f-1}{k}}(\zeta_{p^{f'}-1}^{\frac{p^{f'}-1}{p^f-1}i}+{\mathfrak{P'}})\nonumber\\
&=&\chi_{\mathfrak{P'}}^{p_1}(\zeta_{p^{f'}-1}^{\frac{p^{f'}-1}{kp_1}i}+{\mathfrak{P'}})=\chi_{P'}^{p_1}(\beta^{i}+P'), 
\end{eqnarray}
where $\alpha+P$ and $\beta+P'$ are primitive root of the finite fields
$O_{K}/P$ and $O_{K'}/P'$. 

First, we show the following lemma. 
\begin{lemma}\label{le:Stickel}
Let $\chi_{P'}$ and $\chi_{P}$ be the  Teichm\"{u}ller characters associated to $ P'$ and  $P$, respectively. 
Then, 
\[
(\vartheta_p(\chi_{P'},\chi_P):=)\frac{G_{f'}(\chi_{P'})}{p^{\frac{\phi(k)(p_1-1)}{2e}}G_{f}(\chi_P)}
\] is a $2k'$th or $k'$th root of unity according as 
$k'$ is odd or 
not.   
\end{lemma}
\proof 
First of all, we see 
$\vartheta_p(\chi_{P'},\chi_P)\in \Q(\zeta_{k'})$. 
Note that $\chi_{P}$ is the restriction of $\chi_{P'}$ to $\F_{p^f}$ since 
\[
\chi_{{\mathfrak{P'}}}^{\frac{p^{f'}-1}{k'}}(\zeta_{p^f-1}+{\mathfrak{P'}})=\chi_{{\mathfrak{P}}}^{\frac{p^{f}-1}{k}}(\zeta_{p^f-1}+{\mathfrak{P}})
\]
by $(p^{f'}-1)/k'\equiv (p^{f}-1)/k\,(\mod{p^f-1})$. 
(Thus, in this case, our definition of relative Gauss sums is just the normalization of Yamamoto's relative Gauss sums.)
By Lemma~\ref{Prop:Yamamoto}, $\vartheta_p(\chi_{P'},\chi_P)\in \Q(\zeta_{k'})$ follows. 

Put $f=\frac{\phi(k)}{e}$ and $f'=\frac{\phi(k)p_1}{e}$, and 
set $h=2^t\prod_{i=1}^{\ell}p_i$, where $p_1,p_2,\ldots,p_\ell$ be all distinct prime factors of $k$ and $t$ is the highest power of $2$ dividing $k$. 
It is clear that 
\[
k\sum_{i=0}^{f-1}\left\langle \frac{tp^i}{k}\right\rangle=\sum_{i=0}^{f-1}[tp^i]_k, 
\]
where $[a]_k$ means the reduction of $a$ modulo $k$. 
In other words, it is equal to
\begin{eqnarray*}
\sum_{x\in\langle p\rangle \le (\Z/k\Z)^\ast} [tx]_k&=&\sum_{y=0}^{\frac{k}{h}-1}
\sum_{z\in\langle p\rangle \le (\Z/h\Z)^\ast} hy+ [tz]_{h}\\
&=&k(\frac{k}{h}-1)\phi(h)/2e+\frac{k}{h}
\sum_{z\in\langle p\rangle \le  (\Z/h\Z)^\ast} [tz]_{h}\\
\end{eqnarray*}
where 
note that $t$'s modulo $h$ again forms a set of representatives of $ (\Z/h\Z)^\ast/\langle p\rangle$. Thus, we have 
\[
\sum_{i=0}^{f-1}\left\langle \frac{tp^i}{k}\right\rangle=
\frac{\phi(k)-\phi(h)}{2e}+\frac{1}{h}\sum_{z\in\langle p\rangle \le  (\Z/h\Z)^\ast} [tz]_{h}. 
\]
Similarly, we obtain
\[
\sum_{i=0}^{f'-1}\left\langle \frac{tp^i}{k'}\right\rangle=
\frac{\phi(k')-\phi(h)}{2e}+\frac{1}{h}\sum_{z\in\langle p\rangle \le  (\Z/h\Z)^\ast} [tz]_{h}. 
\]
Hence, by the Stickelberger relation, we obtain 
\[
G_{f'}(\chi_{P'})O_{M'}={\mathfrak{p}'}^{(p-1)((f'-\frac{\phi(k')-\phi(h)}{2e})\sum_{t\in T'}\sigma_t-
\sum_{t\in T'}(\frac{1}{h}\sum_{z\in\langle p\rangle \le (\Z/h\Z)^\ast} [tz]_{h})\sigma_t^{-1})}. 
\]
Furthermore, by noting that 
$\mathfrak{p}O_{M'}={\mathfrak{p}'}$, we have 
\[
G_{f}(\chi_{P})O_{M'}={\mathfrak{p}'}^{(p-1)((f-\frac{\phi(k)-\phi(h)}{2e})\sum_{t\in T'}\sigma_t-
\sum_{t\in T'}(\frac{1}{h}\sum_{z\in\langle p\rangle \le (\Z/h\Z)^\ast} [tz]_{h})\sigma_t^{-1})}. 
\]
Since $pO_{M'}={\mathfrak{p}'}^{(p-1)\sum_{i\in T'}\sigma_t}$, it follows that
$p^{\frac{\phi(k')-\phi(k)}{2e}}G_{f}(\chi_{P})O_{M'}=G_{f'}(\chi_{P'})O_{M'}$, i.e., 
$\vartheta_p(\chi_{P'},\chi_P)$ is a unit of $O_{M'}$. 
But, $\vartheta_p(\chi_{P'},\chi_P)\in \Q(\zeta_{k'})$, and hence it is a unit of 
$O_{K'}$. 
Furthermore, all the conjugates of $\vartheta_p(\chi_{P'},\chi_P)$ in $O_{K'}$
have absolute value $1$. Therefore, $\vartheta_p(\chi_{P'},\chi_P)$ is a root of unity 
in $O_{K'}$, which completes the proof. \qed
\begin{lemma}\label{le:Stickel2}
Let $d=2\gcd{(k',p-1)}$ or $\gcd{(k',p-1)}$ according as $k'$ is odd or even. 
Then, it holds that 
\[
\vartheta_{p}(\chi_{P'},\chi_P)^d=1.
\] 
\end{lemma}
\proof
Define $\sigma\in \Gal(\Q(\zeta_{p},\zeta_{k'})/\Q(\zeta_p))$ by $\sigma(\zeta_{pk'})=\zeta_{pk'}^{k'\ell+p}$, where $\ell$ is the inverse of $k'$ modulo $p$. Let $\psi'$ and $\psi$ be the respective canonical additive characters of $\F_{q'}$ and $\F_q$. 
Then, 
\begin{eqnarray*}
\sigma\left(\frac{G_{f'}(\chi_{P'})}{G_f(\chi_{P})}\right)&=&
\sigma\left(\frac{\sum_{\alpha\in\F_{q'}}\psi'(\alpha)\chi_{P'}(\alpha)}{\sum_{\beta\in\F_{q}}\psi(\beta)\chi_{P}(\beta)}\right)\\
&=&\frac{\sum_{\alpha\in\F_{q'}}\psi'((k'\ell+p)\alpha)\chi_{P'}^{k'\ell+p}(\alpha)}{\sum_{\beta\in\F_{q}}\psi((k'\ell+p)\beta)\chi_{P}^{k'\ell+p}(\beta)}\\
&=&\frac{\sum_{\alpha\in\F_{q'}}\psi'(\alpha)\chi_{P'}^{p}(\alpha)}{\sum_{\beta\in\F_{q}}\psi(\beta)\chi_{P}^{p}(\beta)}\\
&=&\frac{G_{f'}(\chi_{P'}^p)}{G_f(\chi_{P}^p)}=\frac{G_{f'}(\chi_{P'})}{G_f(\chi_{P})}. 
\end{eqnarray*} 
Hence, $\sigma(\vartheta_p(\chi_{P'},\chi_P))=\vartheta_p(\chi_{P'},\chi_P)$. On the other hand, in the case when $k'$ is odd,  
since $\vartheta_{p}(\chi_{P'},\chi_P)^{2}=\zeta_{k'}^s$ for some $s$ by Lemma~\ref{le:Stickel}, 
it follows that $\sigma(\vartheta_p(\chi_{P'},\chi_P)^2)=\vartheta_p(\chi_{P'},\chi_P)^{2(k'\ell+p)}=
\vartheta_p(\chi_{P'},\chi_P)^{2p}$, so $\vartheta_p(\chi_{P'},\chi_P)^{2(p-1)}=1$. Together with 
 $\vartheta_p(\chi_{P'},\chi_P)^{2k'}=1$, we obtain $\vartheta_p(\chi_{P'},\chi_P)^{2\gcd{(k',p-1)}}=1$.  In the case when $k'$ is even,  
since $\vartheta_{p}(\chi_{P'},\chi_P)=\zeta_{k'}^s$ for some $s$ by Lemma~\ref{le:Stickel}, it follows that $\sigma(\vartheta_p(\chi_{P'},\chi_P))=\vartheta_p(\chi_{P'},\chi_P)^{k'\ell+p}=
\vartheta_p(\chi_{P'},\chi_P)^{p}$, so $\vartheta_p(\chi_{P'},\chi_P)^{p-1}=1$. Together with 
 $\vartheta_p(\chi_{P'},\chi_P)^{k'}=1$, we obtain $\vartheta_p(\chi_{P'},\chi_P)^{\gcd{(k',p-1)}}=1$.  
\qed

The following is our main theorem of this section. 
\begin{theorem}\label{thm:Stickel}
If  $k'$ is odd and $\gcd{(k',p-1)}=1$, it holds that $\vartheta_{p}(\chi_{P'},\chi_P)=1$. 
\end{theorem}
\proof 
By Lemma~\ref{le:Stickel2}, we have $\vartheta_{p}(\chi_{P'},\chi_P)=-1$ or $1$. 
We consider the reduction of $p^{\frac{\phi(k')-\phi(k)}{2e}}G_{f}(\chi_{P})\vartheta_p(\chi_{P'},\chi_P)$ modulo 
$\lambda:=1-\zeta_{p_1^{t+1}}$, where $t$ is the highest power of $p_1$ dividing 
$k$. 
It is clear that $p^{\frac{\phi(k)(p_1-1)}{2e}}\equiv 1\,(\mod{\lambda})$. 
Let $h:=k'/p_1^{t+1}$. Since $\chi_P$ and $\chi_{P'}$ can be written as 
$\chi_{P}^{xh}\chi_{P}^{yp_1^{t+1}}$ and 
$\chi_{P'}^{xh}\chi_{P'}^{yp_1^{t+1}}$ for some $x$ and $y$ such that $xh+yp_1^{t+1}\equiv 1\,(\mod{k'})$. 
Then, we have $G_{f'}(\chi_{P'})\equiv G_{f'}(\chi_{P'}^{yp_1^{t+1}})\,(\mod{\lambda})$ and $G_{f}(\chi_{P})\equiv G_{f}(\chi_{P}^{yp_1^{t+1}})\,(\mod{\lambda})$, where both of $\chi_{P'}^{yp_1^{t+1}}$ and $\chi_{P}^{yp_1^{t+1}}$ are 
of order $h$. Now, note that 
\begin{eqnarray*}
\chi_{{\mathfrak{P}}}^{yp_1^{t}\frac{p^{f-1}}{k}}
(\zeta_{p^{fp_1}-1}^{i\frac{p^{fp_1}-1}{p^f-1}}+{\mathfrak{P}})=
\chi_{{\mathfrak{P'}}}^{yp_1^{t}\frac{p^{f-1}}{k}}(\zeta_{p^{fp_1}-1}^{i\frac{p^{fp_1}-1}{p^f-1}}+{\mathfrak{P'}})=
\chi_{{\mathfrak{P'}}}^{yp_1^{t+1}\frac{p^{fp_1}-1}{kp_1}}(\zeta_{p^{fp_1}-1}^{i}+{\mathfrak{P'}}). 
\end{eqnarray*}
By the Davenport-Hasse lifting formula, we have 
\begin{eqnarray*}
G_{f'}(\chi_{P'}^{yp_1^{t+1}})&=&G_{f'}(\chi_{{\mathfrak{P}'}}^{yp_1^{t+1}\frac{p^{fp_1-1}}{kp_1}})\equiv
 (-1)^{p_1-1}(G_{f}(\chi_{{\mathfrak{P}}}^{yp_1^{t}\frac{p^{f-1}}{k}}))^{p_1}\,(\mod{\lambda})\\
&=&(G_{f}(\chi_{P}^{yp_1^{t}}))^{p_1}\equiv 
\chi_P^{-yp_1^{t}}(p_1)G_{f}(\chi_{P}^{yp_1^{t+1}})\,(\mod{\lambda}). 
\end{eqnarray*}
Therefore, by noting that $\chi_P^{-yp_1^{t}}(p_1)=1$, we obtain 
\[
G_f(\chi_{P}^{yp_1^{t+1}})(\vartheta_p(\chi_{P'},\chi_P)-1)\equiv 0\,(\mod{\lambda}). 
\]
If $\vartheta_p(\chi_{P'},\chi_P)=-1$, then $\lambda\,|\,2G_f(\chi_{P}^{yp_1^{t+1}})$. Here, by Lemma~\ref{Prop:Yamamoto}, note
that 
\[
G_f(\chi_{P}^{yp_1^{t+1}})=\sum_{x\in L; \Tr_{p^{f}/p}(x)=1}\chi_{P}^{yp_1^{t+1}}(x)\in \Q(\zeta_{k}), 
\] 
where $L$ is a set of representatives for $\F_{p^{f}}^\ast/\F_{p}^\ast$. 
By taking norms of $\lambda$ and $2G_f(\chi_{P}^{yp_1^{t+1}})$ in $\Q(\zeta_{k'})$, we obtain the contradiction 
that $p_1$ divides $2p$. 
\qed

Next, we treat the case when $2\,||\,k'$ and $\gcd{(k'/2,p-1)}=1$. 
\begin{corollary}\label{cor:Stickel}
Assume that $2\,||\,k',k$ and $\gcd{(k'/2,p-1)}=\gcd{(k/2,p-1)}=1$. 
%
Then, \[\vartheta_p(\chi_{P'},\chi_P)=(-1)^{\frac{(p-1)(p_1-1)\phi(h)}{4e}}, 
\] 
where $h$ is the product of all distinct odd prime factors of $k'$. 
\end{corollary}
\proof
Let $U=\Q(\zeta_{k}^2)$, $U'=\Q(\zeta_{k'}^2)$,  
$\tilde{P}=P\cap O_{U}$, and $\tilde{P'}=P'\cap O_{U'}$. 
Then, $\vartheta_{p}(\chi_{\tilde{P'}},\chi_{\tilde{P}})=1$ by Theorem~\ref{thm:Stickel2}. 
Noting that 
\[
G_f(\chi_{P}^2)=G_f(\chi_{{\mathfrak{P}}}^{2\frac{p^{f}-1}{k}})=
G_f(\chi_{\tilde{P}})
\]
and 
\[
G_{f'}(\chi_{P'}^2)=G_{f'}(\chi_{{\mathfrak{P'}}}^{2\frac{p^{f'}-1}{k'}})=
G_{f'}(\chi_{\tilde{P'}}), 
\]
we have 
\begin{eqnarray*}
\vartheta_p(\chi_{P'},\chi_P)&=&\frac{G_{f'}(\chi_{P'})}{p^{\frac{\phi(k')-\phi(k)}{2e}}G_{f}(\chi_P)}\\
&=&\frac{\chi_{P}^2(2)G_{f'}(\chi_{P'}^2)G_{f'}(\chi_{P'}\eta')G_f(\eta)}{p^{\frac{\phi(k')-\phi(k)}{2e}}\chi_{P'}^2(2)G_f(\chi_P^2)G_f(\chi_P\eta)G_{f'}(\eta')}\\
&=&\frac{\chi_P^2(2)G_{f'}(\chi_{\tilde{P'}})G_{f'}(\chi_{\tilde{P'}}^{(1+k'/2)/2})G_f(\eta)}{p^{\frac{\phi(k')-\phi(k)}{2e}}\chi_{P'}^2(2)G_f(\chi_{\tilde{P}})G_f(\chi_{\tilde{P}}^{(1+k/2)/2})G_{f'}(\eta')}\\
\end{eqnarray*}
where 
$\eta'$ and $\eta$ are the respective quadratic characters of $\F_{q'}$ and 
$\F_q$. Since the restrictions of $\chi_P^2$ and $\chi_{P'}^2$  
to $\F_p^\ast$ 
are trivial, we have $\chi_P^2(2)=\chi_{P'}^2(2)=1$. Furthermore, 
since $\vartheta_{p}(\chi_{\tilde{P'}},\chi_{\tilde{P}})=1$, 
we have $G_{f'}(\chi_{\tilde{P'}})/G_{f}(\chi_{\tilde{P}})=p^{(\phi(k')-\phi(k))/2e}$. 
Now, note that  $(1+k/2)/2=h(s-(p_1-1)k/4h)+[gp^m]_h\in (\Z/\frac{k}{2}\Z)^\ast$ if $(1+k'/2)/2=hs+[gp^m]_h\in (\Z/\frac{k'}{2}\Z)^\ast$ for some $g$ in a set of 
representatives of $(\Z/h\Z)^\ast/\langle p\rangle$ and $0\le s\le k'/2h-1$. Hence, by Theorem~\ref{thm:Stickel} and our assumption, 
we have $G_{f'}(\chi_{\tilde{P'}}^{(1+k'/2)/2})/G_f(\chi_{\tilde{P}}^{(1+k/2)/2})=p^{(\phi(k')-\phi(k))/2e}$. Finally, by
the Davenport-Hasse lifting formula and Lemma~\ref{le:quad}, we have 
\[
\frac{G_{f'}(\eta')}{G_{f}(\eta)}=(-1)^{p_1-1}(G_{f}(\eta))^{p_1-1}=(-1)^{\frac{(p-1)(p_1-1)\phi(h)}{4e}}p^{\frac{\phi(k')-\phi(k)}{2e}},
\]
which shows the assertion. 
\qed

\begin{remark}Let $\epsilon$ denote $(-1)^{\frac{(p-1)(p_1-1)\phi(h)}{4e}}$ or $1$ according as 
$2||k$ and $\gcd{(k'/2,p-1)}=1$ or $2\not|k$ and $\gcd{(k',p-1)}=1$. By Theorem~\ref{thm:Stickel2} and Corollary~\ref{cor:Stickel}, for any $a$ s.t. $\gcd{(a,k')}=1$ it is clear that 
\[
\vartheta_p(\chi_{P'}^a,\chi_{P}^a)=\epsilon
\] 
since $
\sigma(\vartheta_p(\chi_{P'},\chi_{P}))=\vartheta_p(\chi_{P'}^a,\chi_{P}^a)$
and $\sigma(\epsilon)=\epsilon$ for $\sigma\in \Gal(\Q(\zeta_{k'p})/\Q(\zeta_p))$ satisfying  
$\sigma(\zeta_{k'})=\zeta_{k'}^a$. 
\end{remark}
\begin{corollary}\label{cor:last}
Assume that  $k'$ is odd and $\gcd{(k',p-1)}=1$. Then, it holds that  $\vartheta_p(\chi_{P'}^{t},\chi_{P}^{t})=1$ for any 
$t$ such that $p^{s}\not|t$, where $s$ is the highest power of $p_1$ 
dividing $k$.  
\end{corollary}
\proof
Put $t=a\cdot \gcd{(t,k)}$ with $\gcd{(a,k)}=1$. 
Let $r'$ and $r$ be the order of $p$ modulo $k'/\gcd{(t,k')}(=:u')$ and 
modulo  $k/\gcd{(t,k)}(=:u)$. Then, by our assumption, $r'=rp_1$ and $u'=up_1$ follow.
Write $J=\Q(\zeta_{u}),J'=\Q(\zeta_{u'}),H=\Q(\zeta_{p^r-1}),H'=\Q(\zeta_{p^{r'}-1})$
and $R=P\cap O_{J}, R'=P'\cap O_{J'},{\mathfrak{R}}={\mathfrak{P}}\cap O_{H},{\mathfrak{R}}'={\mathfrak{P}}'\cap O_{H'}$.  
Then, we have 
\[
\chi_{{\mathfrak{R}}}^{a\frac{p^{r}-1}{u}}(\zeta_{p^{f}-1}^{i\frac{p^f-1}{p^r-1}}+{\mathfrak{R}})=\chi_{{\mathfrak{P}}}^{a\frac{p^{f}-1}{u}}(\zeta_{p^{f}-1}^i+{\mathfrak{P}})
=\chi_{{\mathfrak{P}}}^{\frac{p^{f}-1}{k}t}(\zeta_{p^{f}-1}^i+{\mathfrak{P}})
\]
and hence $\chi_{{\mathfrak{P}}}^{\frac{p^{f}-1}{k}t}$ is the lift of  $\chi_{{\mathfrak{R}}}^{a\frac{p^{r}-1}{u}}$ to $\F_{p^f}$. 
Similarly, 
$\chi_{{\mathfrak{P}}'}^{\frac{p^{f'}-1}{k'}t}$ is the lift of 
$\chi_{{\mathfrak{R}}'}^{a\frac{p^{r'}-1}{u'}}$ to $\F_{p^{f'}}$. 
Now, by 
the Davenport-Hasse lifting formula, we
have 
\begin{eqnarray*}
\vartheta_p(\chi_{P'}^{t},\chi_{P}^{t})&=&
\frac{G_{f'}(\chi_{P'}^{t})}{p^{\frac{\phi(k')-\phi(k)}{2e}}G_f(\chi_{P}^{t})}
=\frac{G_{f'}(\chi_{{\mathfrak{P}}'}^{\frac{p^{f'}-1}{k'}t})}{p^{\frac{\phi(k')-\phi(k)}{2e}}G_f(\chi_{{\mathfrak{P}}}^{\frac{p^{f}-1}{k}t})}
\\
&=&\frac{(-1)^{f'/r'-1}(G_{r'}(\chi_{{\mathfrak{R}}'}^{a\frac{p^{r'}-1}{u'}}))^{f'/r'}}{(-1)^{f/r-1}p^{\frac{\phi(k')-\phi(k)}{2e}}(G_r(\chi_{{\mathfrak{R}}}^{a\frac{p^{r}-1}{u}}))^{f/r}}\\
&=&
\frac{1}{p^{\frac{\phi(k')-\phi(k)}{2e}}}\cdot \left(\frac{G_{r'}(\chi_{R'}^{a})}{G_r(\chi_{R}^{a})}\right)^{f/r}. 
\end{eqnarray*}
Applying  Theorem~\ref{thm:Stickel}, the above is equal to 
\[
\frac{1}{p^{\frac{\phi(k')-\phi(k)}{2e}}}\cdot \left(p^{\frac{r'-r}{2}} \vartheta_p(\chi_{R'},\chi_R)\right)^{f/r}=1, 
\]
which completes the proof. 
\qed
\section{Constructions of strongly regular graphs and skew Hadamard difference sets}\label{sec:const}
\subsection{General construction}
We first recall the following well-known lemma in the theory of difference sets 
(see e.g., \cite{M94,T65}). 
\begin{lemma}\label{Sec3Le1}
Let $(G, +)$ be an abelian group of odd order $v$, $D$ a subset of $G$ 
of size $\frac{v-1}{2}$.  Assume that $D\cap -D=\emptyset$ and 
$0\not \in D$. Then, $D$ is a skew Hadamard difference set 
in $G$ if and only if 
\[
\chi(D)=\frac{-1\pm \sqrt{-v}}{2}
\] 
for all nontrivial characters $\chi$ of $G$.  On the other hand, assume that $0\not\in D$ and $-D=D$. Then $D$ is a Paley type partial difference set in $G$ if and only if 
\[
\chi(D)=\frac{-1\pm \sqrt{v}}{2}
\] 
for all nontrivial characters $\chi$ of $G$. 
\end{lemma}

Let $q$ be a prime power and let 
$C_i^{(k,q)}=\gamma^i \langle \gamma^k\rangle$, $0\le i\le k-1$, be 
the cyclotomic classes of order $k$ of $\F_q$, where $\gamma$ is a fixed primitive element of $\F_q$. From now on, we will assume that 
$D$ is a union of cyclotomic classes of order $k$ of $\F_q$. In order to check whether a candidate subset, $D=\bigcup_{i\in I}C_i^{(k,q)}$, is a connection set of a strongly regular Cayley graph (i.e., a regular partial difference set), we will compute the sums $\psi(aD):=\sum_{x\in D}\psi(ax)$ for all $a\in \F_q^\ast$,  where $\psi$ is the canonical additive character of $\F_q$,  
since the restricted eigenvalues of Cayley graph $\Cay(\F_q,D)$, as explained in \cite[p.~134]{bh}, are $\psi(\gamma^a D),$ where $a=0,1,\ldots, q-2$. Similarly, to check whether $D$ is a  skew Hadamard difference set in $(\F_q,+)$, we will compute the sums $\psi(aD)$ for all $a\in \F_q^\ast$ because of  Lemma~\ref{Sec3Le1}. Thus, by Theorem~\ref{char} and Lemma~\ref{Sec3Le1}, in both cases
we need to show that the set $\{\psi(\gamma^a D)\,|\,a=0,1,\ldots,q-2\}$ has precisely two elements. 
Note that the sum $\psi(aD)$ can be expressed as a linear combination of 
Gauss sums using the orthogonality of characters: 
\begin{eqnarray*}
\psi(aD)&=&\frac{1}{k}\sum_{i\in I}\sum_{x\in \F_q^\ast}\psi(a\gamma^i x^k)\\
&=&\frac{1}{k}\sum_{i\in I}\sum_{x\in \F_q^\ast}\frac{1}{q-1}
\sum_{y\in \F_q^\ast}\psi(y)
\sum_{\chi\in \widehat{\F_q^\ast}}\chi(a\gamma^ix^k)\overline{\chi(y)}\\
&=&\frac{1}{(q-1)k}\sum_{i\in I}\sum_{x\in \F_q^\ast}
\sum_{\chi\in \widehat{\F_q^\ast}}G(\chi^{-1})
\chi(a\gamma^ix^k)\\
&=&\frac{1}{(q-1)k}\sum_{i\in I}
\sum_{\chi\in \widehat{\F_q^\ast}}G(\chi^{-1})
\chi(a\gamma^i)\sum_{x\in \F_q^\ast}\chi(x^k)\\
&=&\frac{1}{k}
\sum_{\chi\in C_0^{\perp}}G(\chi^{-1})
\sum_{i\in I}\chi(a\gamma^i ), 
\end{eqnarray*}
where $\widehat{\F_q^\ast}$ is the group of multiplicative characters of  
$\F_q^\ast$ and 
$C_0^{\perp}$ is the subgroup of $\widehat{\F_q^\ast}$
consisting of all $\chi$ which are trivial on $C_0^{(k,q)}$. 

In this section, similar to Section~\ref{RGauss}, we 
will assume the following. 
Let $h=2^tp_1p_2\cdots p_\ell$ be a positive integer with distinct odd primes  $p_i$ and let $p$ be a prime satisfying the following:
For any divisor $d=2^s p_{i_1}\cdots p_{i_m}$ of $h$, if 
$\langle p\rangle$ is of index $u$ modulo $d$, then so does $\langle p\rangle$ modulo $d'=2^s p_{i_1}^{x_1}\cdots p_{i_m}^{x_m}$ for any $x_i\ge 1$. 
Let $e$ denotes the index of $\langle p\rangle$ modulo $h$. 
We write $k=\prod_{i=1}^\ell 2^t p_i^{e_i}$ and $k'=kp_1$. 
\begin{theorem}\label{Subsec1Thm1}
Let $q=p^{f}$ and $q'=p^{f'}$, where $f=\phi(k)/e$ and $f'=\phi(k')/e$, and let 
\[
J=\{x\,|\,\mbox{$x$ divides $k$ and  $x$ is not divisible by $p_1^{e_1}$}\}\subseteq \N. 
\]
Let $J_1$ and $J_2$ be a partition of 
$J$ into two parts and let $I$ be a subset of $\{0,1,\ldots,k-1\}$ satisfying the following conditions:
\begin{itemize}
\item[(i)] $\sum_{i\in I}\zeta_{k}^{ij}=0$ for all $j\in J_1$.  
\item[(ii)] 
$\theta_p(\chi_{P'}^{j},\chi_{P}^{j})=\epsilon$ for all 
$j\in J_2$, where $\epsilon=1$ or $-1$ not depending on $j$. 
\item[(iii)] If $\ell\ge 2$ or $t\ge 1$, 
\[
G_{f'}(\chi_{P'}^{-p_1^{e_1+1}v})
=\epsilon p^{\frac{\phi(k')-\phi(k)}{2e}}G_{f}(\chi_{P}^{-p_1^{e_1}v})
\]
for all $1\le v\le k/p_1^{e_1}-1$. 
\end{itemize}
Let 
\[
D=\bigcup_{i\in I}C_i^{(k,q)} \mbox{\, and\, }
D'=\bigcup_{i\in I}\bigcup_{j=0}^{p_1-1}C_{ip_1+jk/p_1^{e_1}}^{(kp_1,q')}.
\] 
Assume that 
the size of the set $\{\psi(\gamma^a D)\,|\,a=0,1,\ldots,q-2\}$ is exactly two, where $\gamma$ is a  primitive root of $\F_q$ and 
$\psi$ is the canonical additive character of $\F_q$.
Then, the size of the set $\{\psi'(\omega^a D')\,|\,a=0,1,\ldots,q'-2\}$  is exactly two, where $\omega$ is a primitive root of $\F_{q'}$ and $\psi'$ is the canonical additive character of $\F_{q'}$. 
\end{theorem}
\proof 
In this proof, without loss of generality, we assume that the primitive roots $\gamma$ and $\omega$ have the forms 
$\gamma=\alpha+P\in O_{K}/P$ and $\omega=\beta+P'\in O_{K'}/P'$ for $\alpha$ and $\beta$ of 
(\ref{eq:basicchar}). Then, $\chi_{P'}^{u}(\omega^{p_1})=\chi_{P}^{u}(\gamma)$ follows. 

To prove the theorem,  
it is sufficient to evaluate the sum
\[
kp_1\cdot \psi'(\omega^a D')=\sum_{u=0}^{kp_1-1}G_{f'}(\chi_{P'}^{-u}) \sum_{i\in I}\sum_{j=0}^{p_1-1}
\chi_{P'}^{u}(\omega^{a+ip_1+jk/p_1^{e_1}}), 
\]
where $a=0,1,\ldots, k'-1$ and $\psi'$ is the canonical additive character of $\F_{q'}$. 

For $u=0$,  
we have 
\begin{eqnarray*}
G_{f'}(\chi_{P'}^0)\sum_{i\in I}\sum_{j=0}^{p_1-1}\chi_{P'}^0(\omega^{a+ip_1+jk/p_1^{e_1}})=-p_1|I|. 
\end{eqnarray*}

For $u=p_1^{e_1}v$ with  $v\not\equiv 0\,(\mod{p_1})$,  
we have 
\begin{equation}\label{eq:gene3}
G_{f'}(\chi_{P'}^{-p_1^{e_1}v}) \sum_{i\in I}\sum_{j=0}^{p_1-1}
\chi_{P'}^{p_1^{e_1}v}(\omega^{a+ip_1+jk/p_1^{e_1}})=0. 
\end{equation}

If $\ell\ge 2$ or $t\ge 1$, for $u=p_1^{e_1+1}v$ with  $v\not=0$,  
we have 
\begin{eqnarray*}
G_{f'}(\chi_{P'}^{-p_1^{e_1+1}v}) \sum_{i\in I}\sum_{j=0}^{p_1-1}
\chi_{P'}^{p_1^{e_1+1}v}(\omega^{a+ip_1+jk/p_1^{e_1}})
=p_1G_{f'}(\chi_{P'}^{-p_1^{e_1+1}v})\sum_{i\in I}
\chi_{P'}^{p_1^{e_1+1}v}(\omega^{a+ip_1+jk/p_1^{e_1}})
\end{eqnarray*}
for any $j$. 
Note that for each $a\in \{0,1,\ldots,k'-1\}$, there is a unique 
$j\in \{0,1,\ldots,p_1-1\}$ such that 
$p_1\,|\,a+jk/p_1^{e_1}$; we write $a+jk/p_1^{e_1}=p_1j_a$. 
Then, the above is equal to 
\begin{equation}\label{eq:gene2}
p_1G_{f'}(\chi_{P'}^{-p_1^{e_1+1}v})\sum_{i\in I}
\chi_{P'}^{p_1^{e_1+1}v}(\omega^{p_1(j_a+i)}). 
\end{equation}
Furthermore, 
since 
$\chi_{P'}^{u}(\omega^{p_1(j_a+i)})=\chi_{P}^{u}(\gamma^{j_a+i})$,
by the assumption (iii), eq.~(\ref{eq:gene2}) is rewritten as 
\begin{equation}\label{eq:gene4}
\epsilon p_1p^{\frac{\phi(k')-\phi(k)}{2e}}G_{f}(\chi_{P}^{-p_1^{e_1}v})\sum_{i\in I}
\chi_{P}^{p_1^{e_1}v}(\gamma^{j_a+i}). 
\end{equation}

For the remaining cases, we can assume that $p_1^{e_1}\not | u$, and write 
$u=kv_1+v_2$ for some  $0\le v_1\le p_1-1$ and  $0\le v_2\le k-1$. 
Then, since $G_{f'}(\chi_{P'}^{kv_1+v_2}) =G_{f'}(\chi_{P'}^{kv_1'+v_2}) $ for 
$0\le v_1,v_1'\le p_1-1$, 
we have 
\begin{eqnarray*}
& &\sum_{v_1=0}^{p_1-1}\sum_{v_2=1}^{k-1}
G_{f'}(\chi_{P'}^{-kv_1-v_2}) \sum_{i\in I}\sum_{j=0}^{p_1-1}
\chi_{P'}^{kv_1+v_2}(\omega^{a+ip_1+jk/p_1^{e_1}})\\&=&
p_1\sum_{v_2=1}^{k-1}
G_{f'}(\chi_{P'}^{-v_2})\sum_{i\in I}
\chi_{P'}^{v_2}(\omega^{p_1(j_a+i)})\\
&=&p_1\sum_{v_2=1;\gcd{(v_2,k)}\in J_2}^{k-1}
G_{f'}(\chi_{P'}^{-v_2})\sum_{i\in I}
\chi_{P'}^{v_2}(\omega^{p_1(j_a+i)}). 
\end{eqnarray*}
By our assumption that $G_{f'}(\chi_{P'}^{-v_2})=\epsilon p^{\frac{\phi(k')-\phi(k)}{2e}}G_{f}(\chi_{P}^{-v_2})$ and by 
$\chi_{P'}^{v_2}(\omega^{p_1(j_a+i)})=\chi_{P}^{v_2}(\gamma^{j_a+i})$,
the above is 
equal to  
\begin{eqnarray*}
\epsilon p_1p^{\frac{\phi(k')-\phi(k)}{2e}}\sum_{v_2=1;\gcd{(v_2,k)}\in J_2}^{k-1}G_{f}(\chi_{P}^{-v_2})
\sum_{i\in I}
\chi_{P}^{v_2}(\gamma^{j_a+i}). 
\end{eqnarray*}
Finally, together with eq.~(\ref{eq:gene3}) and (\ref{eq:gene4}), 
we obtain 
\[
kp_1\cdot \psi'(\omega^a D')+p_1|I|=\epsilon p_1p^{\frac{\phi(k')-\phi(k)}{2e}}\sum_{v_2=1}^{k-1}G_{f}(\chi_{P}^{-\ell_2})
\sum_{i\in I}
\chi_{P}^{\ell_2}(\gamma^{j_a+i}). 
\]
Now, by the assumption that the size of the set 
\[
\{\psi(\gamma^a D)\,|\,a=0,1,\ldots,k-1\}
\]
is exactly two, we obtain the assertion. In particular, the two values 
in $\{\psi'(\omega^a D')\,|\,a=0,1,\ldots,q'-2\}$ are given as 
\begin{equation}\label{eq:gene}
\frac{1}{kp_1}(\epsilon p_1p^{\frac{\phi(k)(p_1-1)}{2e}}(ks+|I|)-p_1|I|)=
\epsilon p^{\frac{\phi(k)(p_1-1)}{2e}}s+\frac{|I|(\epsilon p^{\frac{\phi(k)(p_1-1)}{2e}}-1)}{k}, 
\end{equation}
where $s=\psi(\gamma^a D)$ for some $a$.  
\qed
\subsection{Strongly regular graphs}\label{sec:stro}
In this subsection, we write $k=\prod_{i=1}^\ell p_i^{e_i}$, where $p_i$ are distinct odd primes and  assume that $p$ is a prime such that  $\ord_{k}(p)=\phi(k)/e$. Furthermore, assume that 
$\langle p\rangle$ is again of index $e$ modulo $k'(:=kp_1)$ and 
$\gcd{(k',p-1)}=1$. 
\begin{theorem}\label{cor:srg1}
Let $h=p_1\cdots p_mp_{m+1}\cdots p_\ell$ with all distinct odd primes $p_i$ and $[(\Z/h\Z)^\ast:\langle p\rangle]=e$.  Furthermore, 
Let $k=p_1^{e_1}\cdots p_m^{e_m} p_{m+1}^{e_{m+1}}\cdots p_\ell^{e_{\ell}}$, where 
$e_i\ge 1$ for $1\le i\le m$ and $e_i= 1$ for $m+1\le i\le \ell$, and assume that $\langle p\rangle $ is again of index $e$ modulo $k$. 
Let $q_1=p^{d}$ and $q=p^{f}$, where $d=\phi(h)/e$ and  $f=\phi(k)/e$. 
Put
$h_j=\prod_{i\not=j} p_i$ for $1\le j\le m$. 
Assume that there exists an integer $s_j$ s.t. 
$p^{s_j}\equiv -1\,(\mod{h_j})$ for $1\le j\le m$. 
Let 
\[
D:=\bigcup_{i_1=0}^{p_1^{e_1-1}-1}\cdots \bigcup_{i_m=0}^{p_m^{e_m-1}-1}C_{i_1{n}_1+\cdots+i_m {n}_m }^{(k,q)}, 
\]  
where $n_j=\prod_{i\not=j} p_i^{e_i}$. 
If $\Cay(\F_{q_1},C_0^{(h,q_1)})$ is an srg, then  so 
does $\Cay(\F_q,D)$.
\end{theorem}
\proof
We will show by induction. 
Write  
\[
D=\bigcup_{i_1=0}^{p_1^{e_1-1}-1}\cdots \bigcup_{i_m=0}^{p_m^{e_m-1}-1}C_{i_1{n}_1+\cdots+i_m {n}_m }^{(k,q)} 
\] 
and 
assume the size of the set $\{\psi(\gamma^a D)\,|\,a=0,1,\ldots,q-2\}$  is exactly two. We put 
\[
I=\bigcup_{i_1=0}^{p_1^{e_1-1}-1}\cdots \bigcup_{i_m=0}^{p_m^{e_m-1}-1}\{i_1{n}_1+\cdots+i_m {n}_m\}
\]
in Theorem~\ref{Subsec1Thm1}. Let
$J$ be the set of positive divisors of $k$ not divisible by $p_1^{e_1}$, 
\[
J_1=\{x\,|\,\mbox{$\exists i$, $1\le i\le m$, s.t. $p_i^r\,||\,x$, where $1\le r\le e_i-1$}\}\subseteq J, 
\]
and $J_2=J\setminus J_1$. 
Then, by the definition of $I$, 
it is clear that $\sum_{i\in I}\zeta_{k}^{ij}=0$ for all $j\in J_1$. 

Furthermore, since the assumption $p^{s_1}\equiv -1\,(\mod{h_1})$ implies that  $p$ is semi-primitive modulo $n_1$, by Theorems~\ref{thm:semiprim}, for $u=p_1^{e_1+1}v$ we have 
\[
G_{f'}(\chi_{P'}^{p_1^{e_1+1}v})=p^{\frac{\phi(k')-\phi(k)}{2e}}
G_{f}(\chi_{P}^{p_1^{e_1}v}). 
\]

Moreover, by Corollary~\ref{cor:last}, we have for any $a\in J_2$ 
\[
G_{f'}(\chi_{P'}^a)=p^{\frac{\phi(k')-\phi(k)}{2e}}G_{f}(\chi_{P}^a).
\]
Thus, the assumptions (i), (ii), and (iii) of Theorem~\ref{Subsec1Thm1} are satisfied. 
Now, by applying Theorem~\ref{Subsec1Thm1},  the size of the set $\{\psi'(\gamma^a D')\,|\,a=0,1,\ldots,q'-2\}$  is exactly two, 
where 
\begin{eqnarray*}
D'&=&\bigcup_{j=0}^{p_1-1}\bigcup_{i_1=0}^{p_1^{e_1-1}-1}\cdots \bigcup_{i_m=0}^{p_m^{e_m-1}-1}C_{p_1(i_1{n}_1+\cdots+i_m {n}_m )+jn_1}^{(kp_1,q')}\\
&=&\bigcup_{i=0}^{p_1^{e_1}-1}\bigcup_{i_2=0}^{p_2^{e_2-1}-1}\cdots \bigcup_{i_m=0}^{p_m^{e_m-1}-1}C_{in_1+i_2{n'}_2+\cdots+i_m {n'}_m }^{(kp_1,q')}
\end{eqnarray*}
with  ${n'}_i=n_ip_1$. 
\qed
\begin{example}
\begin{itemize}
\item[(i)] If $\ell=1$ in Theorem~\ref{cor:srg1}, we do not need the condition that 
there exists an integer $s_j$ s.t. 
$p^{s_j}\equiv -1\,(\mod{h_j})$. Hence, assuming that 
\[
 [(\Z/p_1\Z)^\ast:\langle p\rangle]=[(\Z/p_1^{e_1} \Z)^\ast:\langle p\rangle]=e, 
\]
if 
$\Cay(\F_{p^{\phi(p_1)/e}},C_0^{(p_1,p^{\phi(p_1)/e})})$ forms an srg, 
then so does $\Cay(\F_{p^{\phi(p_1^{e_1})/e}},D)$, where 
\[
D=\bigcup_{i=0}^{p_1^{e_1-1}-1}C_i^{(p_1^{e_1},p^{{\phi(p_1^{e_1})/e}})}.
\] It is easy to see by induction that 
$\ord_{p_1^{e_1}}(p)=\phi(p_1^{e_1})/e$ for general $e$ and for all pairs $(k=p_1,p)$ of No. 1, 2, 4, 5, 6, 7, 9, and 11 in Table~\ref{Tab1}. Thus, all these srgs can be generalized into infinite families. Note that there are a lot of examples in subfield case satisfying $[(\Z/p_1\Z)^\ast:\langle p\rangle]=e$ and 
$p_1=\frac{p^{\phi(p_1)/e}-1}{p^t-1}$ for some $t\,|\,\phi(p_1)/e$. For example, we list ten examples
satisfying these conditions in Table~\ref{Tab4}. 
\begin{table}[h]
\caption{Subfield examples of $\ell=1$ led to infinite families  }
\label{Tab4}
$$
\begin{array}{|c|c|c|c|}
\hline
p_1&p&f&e:=[(\Z/k\Z)^\ast:\langle p\rangle]\\
\hline
7&2&2&2\\
13&3&3&4\\
31&2&5&6\\
31&5&3&10\\
73&2&9&8\\
127&2&7&18\\
307&17&3&102\\
757&3&9&84\\
1093&3&7&156\\
1723&41&3&574\\
\hline
\end{array}
$$
\end{table}
These examples can be similarly generalized into nontrivial infinite families. 
\item[(ii)] If $\ell=2$, 
 in Theorem~\ref{cor:srg1}, we need the condition that 
there exists an integer $s_i$ s.t. 
$p^{s_i}\equiv -1\,(\mod{p_i})$ for either of $i=1,2$. 
Hence, assuming that  
$p$ is semi-primitive modulo both of $p_1$ and $p_2$,  and 
\[
[(\Z/p_1p_2\Z)^\ast:\langle p\rangle]=[(\Z/p_1^{e_1}p_2^{e_2} \Z)^\ast:\langle p\rangle]=e,
\]
if  
$\Cay(\F_{p^{\phi(p_1p_2)/e}},C_0^{(p_1p_2,p^{\phi(p_1p_2)/e})})$ forms an srg, 
then so does $\Cay(\F_{p^{\phi(p_1^{e_1}p_2^{e_2})/e}},D)$, where \[
D=\bigcup_{i=0}^{p_1^{e_1-1}-1}\bigcup_{j=0}^{p_2^{e_2-1}-1}
C_{i_1p_2^{e_2}+i_2 p_1^{e_1} }^{(p_1^{e_1}p_2^{e_2},p^{\phi(p_1^{e_1}p_2^{e_2})/e})}. 
\] 
It is easy to see by induction that 
$\ord_{p_1^{e_1}p_2^{e_2}}(p)=\phi(p_1^{e_1}p_2^{e_2})/e$ for any $e_1,e_2$ and for pairs $(k=p_1p_2,p)$ of No. 3 and 10 in Table~\ref{Tab1}. Thus, these srgs can be generalized into infinite families. 
On the other hand, 
if $p$ is semi-primitive modulo either one of $p_1$ or $p_2$, say $p_2$, 
then $\Cay(\F_{p^{\phi(p_1^{e_1}p_2)/e}},D)$  forms an srg
under the assumption that 
$[(\Z/p_1^{e_1}p_2 \Z)^\ast:\langle p\rangle]=e$, where \[
D=\bigcup_{i=0}^{p_1^{e_1-1}-1}
C_{i_1p_2}^{(p_1^{e_1}p_2,p^{\phi(p_1^{e_1}p_2)/e})}. 
\] 
It is easy to see by induction that 
$\ord_{p_1^{e_1}p_2}(p)=\phi(p_1^{e_1}p_2)/e$ for any $e_1$ and $p$ is semi-primitive modulo $p_2$ for the triple $(p_1,p_2,p)=(19,7,5)$ of No. 8 in Table~\ref{Tab1}. Thus, this srg can be generalized into infinite families. 
Moreover, we can find some examples in subfield case satisfying $[(\Z/p_1p_2\Z)^\ast:\langle p\rangle]=e$ and 
$p_1p_2=\frac{p^{\phi(p_1p_2)/e}-1}{p^t-1}$ for some $t\,|\,\phi(p_1p_2)/e$. For example, we list 
four examples
satisfying these conditions in Table~\ref{Tab6}. 
\begin{table}[h]
\caption{Subfield examples of $\ell=2$ led to infinite families}
\label{Tab6}
$$
\begin{array}{|c|c|c|c|c|c|}
\hline
p_1&p_2&p&f&e:=[(\Z/k\Z)^\ast:\langle p\rangle]&\mbox{sp}\\
\hline
3&5&2&4&2&\mbox{b}\\
5&17&2&8&8&\mbox{b}\\
31&11&2&10&30&\mbox{o}\\
127&43&2&14&378&\mbox{o}\\
\hline
\end{array}
$$
\end{table}
In the sixth column ``sp''  of the table,  ``b'' indicates that 
$p$ is semi-primitive modulo both of $p_1$ and $p_2$, and 
``o'' indicates that $p$ is semi-primitive modulo $p_2$ only. 
These examples can be generalized into nontrivial infinite families. 
\end{itemize}
\end{example}
\subsection{Skew Hadamard difference sets}
In this subsection, we write $k=2 p_1^{e_1}$, where $p_1$ is an odd prime and  assume that $p$ is a prime such that  $\ord_{k}(p)=\phi(k)/e$. Furthermore, assume that 
$p$ is again of index $e$ modulo $k'(:=kp_1)$ and $\gcd{(k'/2,p-1)}=1$. 
\begin{theorem}\label{cor:sHd1}
Let $h=2p_1$ with an odd prime $p_1$ and let $p$ be a prime such that $\langle p\rangle$ is of index $e$ modulo $h$. Furthermore, 
let $k=2p_1^{e_1}$ and assume that $\langle p\rangle$ is again of index $e$ modulo $k$. 
Put $q_1=p^{d}$ and $q=p^{f}$, where $d=\phi(h)/e$ and  $f=\phi(k)/e$. 
Define $H$ as any subset of $\{0,1,\ldots,h-1\}$ such that $\sum_{i\in H}\zeta_{p_1}^{i}=0$.
Let \[
D=\bigcup_{i\in H}C_{i}^{(h,q_1)} \mbox{\, and \, } 
D'=\bigcup_{i_1=0}^{p_1^{e_1-1}}\bigcup_{i\in H}C_{2i_1+ik/h}^{(k,q)}. 
\] 
If $D$ is a skew Hadamard difference set or a Paley type regular partial difference set  on $\F_{q_1}$, then  so does 
$D'$  on $\F_q$.
\end{theorem}
\proof
We will show by induction. 
Write 
\[
D=\bigcup_{i_1=0}^{p_1^{e_1-1}}\bigcup_{i\in H}C_{2i_1+ik/h}^{(k,q)}
\]
and assume that 
the size of the set $\{\psi(\gamma^a D)\,|\,a=0,1,\ldots,q-2\}$  is exactly two, which are $
\frac{-1\pm \sqrt{\tau q}}{2}$, 
where $\tau=1$ 
or $-1$ according as $D$ is a Paley type regular partial difference set or a skew Hadamard difference set. 
Now, we  put 
\[
I=\bigcup_{i_1=0}^{p_1^{e_1-1}-1}\bigcup_{i\in H}
\{2i_1+ik/h\}
\]
in Theorem~\ref{Subsec1Thm1}. 
Let 
\[
J_2=\{1\} \subseteq J=\{1,p_1,\ldots,p_1^{e-1}\}\cup 2\{1,p_1,\ldots,p_1^{e-1}\}
\]
and $J_1=J\setminus J_2$. 
Then, by the definition of $I$, 
it is clear that $\sum_{i\in I}\zeta_{k}^{ij}=0$ for all $j\in J_1$. 
By Lemma~\ref{le:quad}, we have 
\[
G_{f'}(\chi_{P'}^{k'/2})=(-1)^{\frac{(p-1)(p_1-1)\phi(h)}{4e}}p^{\frac{\phi(k')-\phi(k)}{2e}}G_{f}(\chi_{P}^{k/2}). 
\]
Furthermore, by Corollary~\ref{cor:Stickel}, we have 
\[
G_{f'}(\chi_{P'})=(-1)^{\frac{(p-1)(p_1-1)\phi(h)}{4e}}p^{\frac{\phi(k')-\phi(k)}{2e}}G_{f}(\chi_{P}). 
\]
Thus, the assumptions (i),  (ii), and (iii) of Theorem~\ref{Subsec1Thm1} are satisfied. 
Now, by applying Theorem~\ref{Subsec1Thm1},  the size of the set $\{\psi'(\gamma^a D')\,|\,a=0,1,\ldots,q'-2\}$  is exactly two, 
where 
\begin{eqnarray*}
D'&=&\bigcup_{j=0}^{p_1-1}\bigcup_{i_1=0}^{p_1^{e_1-1}}\bigcup_{i\in H}C_{(2i_1+ik/h)p_1+jk/p_1^{e_1}}^{(k,q)}\\
&=&\bigcup_{i_1=0}^{p_1^{e_1}-1}\bigcup_{i\in H}C_{2i_1+ik'/h}^{(k,q)}.
\end{eqnarray*}
In particular, by eq.~(\ref{eq:gene}), the two values in $\{\psi'(\gamma^a D')\,|\,a=0,1,\ldots,q'-2\}$ are 
\[
\epsilon p^{\frac{\phi(k)(p_1-1)}{2e}}\left(\frac{-1\pm \sqrt{\tau p^f}}{2}\right)+\frac{k}{2}\cdot \left(\frac{\epsilon p^{\frac{\phi(k)(p_1-1)}{2e}}-1}{k}\right)=\frac{-1\pm \epsilon \sqrt{\tau p^{f'}}}{2}, 
\]
which completes the proof. 
\qed
\begin{example}
In \cite{FMX11}, several examples satisfying the condition of  
Theorem~\ref{cor:sHd1} were found from index $2$ case, which were 
generalized into infinite families using Gauss sums of index $2$. 
We can find by computer further two examples having the following parameters from index $4$ case: 
\[
(p_1,p,f,e)=(13,3,3,4) \mbox{ and } (29,7,7,4). 
\] 
In particular, the latter example was found by Tao Feng \cite{F12}. 
We choose $H$ in Theorem~\ref{cor:sHd1} as 
$H=Q \cup 2Q\cup \{p_1\}$ for the 
former parameter and choose
$H=Q \cup 2Q\cup \{0\}$ for the latter 
parameter, where $Q$ is the subgroup of index $2$ of $(\Z/2p_1\Z)^\ast$. It is easy to check that these $H$ satisfies the condition of 
Theorem~\ref{cor:sHd1} and $\langle p\rangle$ is of index $4$ in 
$(\Z/2p_1^{e_1}\Z)^\ast$ for general $e_1$. Hence, 
these examples can be generalized into infinite families. 
\end{example}
\section{Final remarks}

We close this paper by referring the reader to the interesting paper \cite{Wu12} by Wu. Immediately after writing up this manuscript, the author 
became aware that Wu \cite{Wu12} obtained a nice result on the existence problem of cyclotomic srgs. 

In our paper,  cyclotomic constructions of strongly regular Cayley graphs and skew Hadamard difference sets on $\F_q$ were given. For example, we proved the following result (which follows from the more general theorem \ref{cor:srg1}): 
For an odd prime $p_1$, assume that (i) $\gcd{(p(p-1),p_1)}=1$ (ii) $\langle p\rangle$ is 
of index $e$ modulo $p_1$ (iii) $\Cay(\F_{p^{(p_1-1)/e}},C_0^{(p^{(p_1-1)/e},p_1)})$ is strongly regular. 
Then, if $\langle p\rangle$ is of index $e$ modulo $p_1^m$,
$\Gamma=\Cay(\F_{p^{p_1^{m-1}(p_1-1)/e}},\bigcup_{i=0}^{p_1^{m-1}-1}C_i^{(p^{p_1^{m-1}(p_1-1)/e},p_1^m)})$ 
is also strongly regular.  Since there are a lot of subfield or sporadic examples satisfying the assumption of this result, we 
consequently obtain many new infinite families of strongly regular Cayley graphs. 
This result can be viewed as a ``recursive'' construction of srgs not saying
anything about the existence of ``starting'' srgs. 

On the other hand,  Wu \cite{Wu12} gave necessary and sufficient conditions for 
$\Gamma$ to be an srg by generalizing the method used in the paper of Ge, Xiang, and Yuan \cite{GXY11}. 
Although it seems that the assumptions of our main result are simpler and the situation is definitely much more general than that of \cite{Wu12}, the approach in \cite{Wu12} is obviously different from ours and  his results are not completely included in ours.  In fact, Wu \cite{Wu12} obtained two conditions (one is an equation and  the other is a congruence) which are necessary and sufficient
for the construction to give rise to an srg, and  his approach has the advantage of revealing an interesting 
connection between strongly regular Cayley graphs $\Cay(\F_{p^{(p_1-1)/e}},C_0^{(p^{(p_1-1)/e},p_1)})$ and 
cyclic difference sets in $(\Z/p_1\Z, +)$, which will be very effective to get some new cyclic difference 
sets and also a strong necessary condition for the existence of cyclotomic srgs.  
\section*{Acknowledgements} 
The work of K. Momihara was supported by JSPS under Grant-in-Aid for Research Activity Start-up 23840032. 

\end{document}